\documentclass[11pt,english,reqno]{amsart}

\usepackage{amsthm,amsmath,amssymb,amsfonts}
\usepackage{mathtools,mleftright}
\usepackage{microtype,setspace,epigraph}
\usepackage[dvipsnames]{xcolor}
\usepackage{graphicx,tikz}
\usetikzlibrary{arrows.meta,positioning,calc}
\usepackage{booktabs,float,caption,subcaption} 
\usepackage[a4paper,margin=1in]{geometry}

\definecolor{webgreen}{rgb}{0,.5,0}
\definecolor{webbrown}{rgb}{.6,0,0}

\usepackage[
  pdfauthor={},
  pdfkeywords={},
  pdftitle={},
  pdfcreator={},
  pdfproducer={},
  colorlinks,linkcolor=blue,filecolor=webbrown,citecolor=webgreen]{hyperref}

\usepackage{cleveref}
\mleftright

\theoremstyle{plain}
\newtheorem{theorem}{Theorem}[section]
\newtheorem{proposition}[theorem]{Proposition}
\newtheorem{lemma}[theorem]{Lemma}
\newtheorem{corollary}[theorem]{Corollary}

\theoremstyle{definition}
\newtheorem{example}[theorem]{Example}
\newtheorem{remark}[theorem]{Remark}

\numberwithin{equation}{section}
\numberwithin{table}{section}
\numberwithin{figure}{section}

\newcommand{\ZZ}{\mathbb{Z}}
\newcommand{\QQ}{\mathbb{Q}}

\newcommand{\CC}{\mathbb{C}}

\newcommand{\AP}{\mathbf{A}}
\newcommand{\AR}{\mathbf{A}^{\!\!\star}}
\newcommand{\BN}{\mathbf{B}}
\newcommand{\BR}{\mathbf{B}^{\star}}
\newcommand{\EN}{\mathbf{E}}
\newcommand{\ER}{\mathbf{E}^{\star}}
\newcommand{\FN}{\mathcal{F}}
\newcommand{\GN}{\mathbf{G}}
\newcommand{\GD}{\widehat{\mathbf{G}}}

\newcommand{\RT}{\mathfrak{R}}
\newcommand{\CT}{\mathfrak{C}}
\newcommand{\ST}{\mathfrak{S}}

\newcommand{\mids}{\,\mid\,}
\newcommand{\nmids}{\,\nmid\,}
\newcommand{\textq}[1]{\quad \text{#1} \quad}
\newcommand{\andq}{\textq{and}}

\newcommand{\set}[1]{\left\{#1\right\}}
\newcommand{\valueat}[1]{{\;}{\big|{}_{#1}}}
\newcommand{\ditto}{\texttt{"}}

\newcommand{\seqnum}[1]{\href{https://oeis.org/#1}{\rm \underline{#1}}}

\newenvironment{proofcref}[1]{%
  \begin{proof}[Proof of \Cref{#1}]%
}{%
  \end{proof}%
}

\begin{document}

\title[From hockey stick to boomerang via Appell polynomials]
{From hockey stick to boomerang\\ via Appell polynomials}
\author{Bernd C. Kellner}
\address{G\"ottingen, Germany}
\email{bk@bernoulli.org}

\begin{abstract}
In Pascal's triangle, the hockey-stick identities are well known as certain sums 
over the binomial coefficients, e.g., lying in a fixed row or column of the triangle.
We establish unified identities that are connected with Appell polynomials, 
from which the hockey-stick identities follow as a very special case.
Subsequently, we derive several identities for (reciprocal) Appell polynomials 
yielding applications to the Bernoulli and Euler polynomials. As a result, 
we obtain new identities for the Bernoulli, Euler, Fibonacci, and Genocchi numbers, 
as well as, power sums. This also leads to infinitely many recurrence relations 
for these numbers. 
\end{abstract}

\maketitle

% ==============================================================================
\section{Introduction}
% ==============================================================================

Around 1654, Pascal~\cite{Pascal:1665} (posthumously published in 1665) gave a
summation formula of binomial coefficients over the rows that
\begin{equation} \label{eq:hsl}
  \sum_{\nu=0}^{k} \binom{n-k+\nu}{\nu} = \binom{n+1}{k}
\end{equation}
for ${n \geq k \geq 0}$. This formula was independently found by Fermat 
(cf.~Hofmann \cite[p.\,32]{Hofmann:1957}), 
who was in correspondence with Pascal at that time. By symmetry of the binomial 
coefficients, one obtains a sum over the rows but in a fixed column that
\begin{equation} \label{eq:hsr}
  \sum_{\nu=k}^{n} \binom{\nu}{k} = \binom{n+1}{k+1}.
\end{equation}

Today, these sums are known as the \emph{hockey-stick identities}.
In \emph{Pascal's triangle} (see Edwards~\cite{Edwards:2019} for historical notes 
of \cite{Pascal:1665}), the corresponding shapes look like hockey sticks, 
see Figure~\ref{fig:hs1} below. Note that it is convenient to use $n+1$ instead 
of $n$ for the last row to simplify formulas and indices. 

\begin{figure}[H]
\begin{center}
\includegraphics[width=5cm]{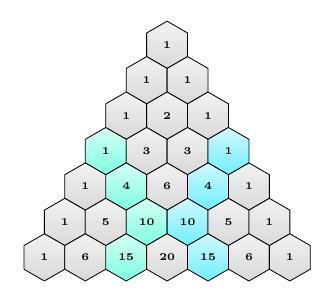}
\caption{Symmetric shapes of hockey sticks in Pascal's triangle. 
The colored shapes on the left and right correspond to 
\eqref{eq:hsl} and \eqref{eq:hsr}, respectively.} 
\label{fig:hs1}
\end{center}
\end{figure}

\begin{figure}[H]
\begin{center}
\includegraphics[width=5cm]{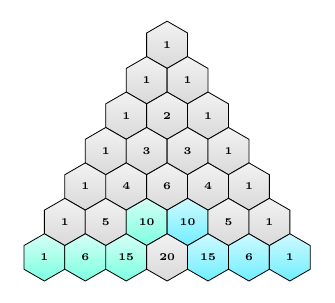}
\caption{Symmetric shapes of hockey sticks in Pascal's triangle. 
The colored shapes on the left and right correspond to 
\eqref{eq:asl} and \eqref{eq:asr}, respectively.} 
\label{fig:hs2}
\end{center}
\end{figure}

In a similar way, an alternating sum in a fixed row is given by
\begin{align}
  \sum_{\nu=0}^{k} \binom{n+1}{\nu} (-1)^{k-\nu} &= \binom{n}{k}. 
  \label{eq:asl} 
\shortintertext{From the symmetry, it follows that}
  \sum_{\nu=k}^{n} \binom{n+1}{\nu+1} (-1)^{\nu-k} &= \binom{n}{k}. 
  \label{eq:asr}
\end{align}
These alternating sums can be also identified by (lying) hockey sticks, 
see Figure~\ref{fig:hs2}.

The sums \eqref{eq:hsr} and \eqref{eq:asl} are given in Comtet 
\cite[Sec.\,1.5, Thm.~E, p.\,10]{Comtet:1974} as two of four basic 
recurrence relations of the binomial coefficients. The remaining relations are
\begin{align}
  \binom{n+1}{k+1} &= \binom{n}{k} + \binom{n}{k+1} 
  \label{eq:bin-add} 
\shortintertext{and}
  \binom{n+1}{k+1} &= \frac{n+1}{k+1} \binom{n}{k}. 
  \nonumber 
\shortintertext{Implicitly, we have used the symmetry rule}
  \binom{n}{k} &= \binom{n}{n-k}  
  \label{eq:bin-sym}
\end{align}
above. As usual, define ${\binom{n}{k} = 0}$ for ${k < 0}$ and ${k > n}$.
See also Graham et al.~\cite[Sec.\,5.1, pp.\,153--165]{GKP:1994}.

To formulate the main results in the next section, we need the following as the 
key ingredient. A sequence $(\AP_n(x))_{n \geq 0}$ is called an Appell sequence 
(named after Appell~\cite{Appell:1880}), where its members are the 
Appell polynomials satisfying the equivalent properties
\begin{align}
  \AP{\!'}_n(x) &= 
  \begin{cases}
    0, & \text{if $n=0$}, \\
    n \AP_{n-1}(x), & \text{otherwise}, 
  \end{cases}
  \label{eq:apd} 
\shortintertext{and}
  \AP_n(x+y) &= \sum_{\nu=0}^{n} \binom{n}{\nu} \AP_{n-\nu}(x) \, y^\nu. 
  \label{eq:apt}
\end{align}
A certain class of the Appell polynomials also obeys the reflection relation
\begin{equation} \label{eq:apr}
  \AP_n(1-x) = (-1)^n \AP_n(x). 
\end{equation}

The prototype of the Appell polynomials is given by the sequence of the monomials 
$(x^n)_{n \geq 0}$. Further important examples are the Bernoulli and Euler polynomials 
$\BN_n(x)$ and $\EN_n(x)$, respectively. These polynomials are introduced later, 
which also obey the reflection relation~\eqref{eq:apr}.

The rest of the paper is organized as follows. The next two sections present the 
main results and several applications. Section~\ref{sec:prel} gives some preliminaries,
while Section~\ref{sec:proofs} contains the proofs. Subsequently, the final section 
concludes the paper with a few remarks.

% ==============================================================================
\section{Main results}
% ==============================================================================

The two different kinds of the hockey-stick identities, as shown in Figures~\ref{fig:hs1} 
and~\ref{fig:hs2}, can be combined to unified identities, respectively. Moreover,
the following theorem connects these identities with Appell polynomials.
The corresponding shapes of the occurring  binomial coefficients in Pascal's triangle 
then look like \emph{boomerangs}, as displayed in Figure~\ref{fig:hs3} below.
In the following, let $x,y \in \CC$. We use the usual convention for power series 
that $x^0 = 1$ for all $x \in \CC$, so in particular $0^0 = 1$.

\begin{theorem} \label{thm:main}
Let $(\AP_\nu(x))_{\nu \geq 0}$ be an Appell sequence.
For $n \geq k \geq 0$, we have
\begin{align}
  \sum_{\nu=0}^{k} \binom{n+1}{\nu} \AP_{k-\nu}(x) \, y^{\nu}
  &= \sum_{\nu=0}^{k} \binom{n-\nu}{k-\nu} \AP_\nu(x+y) \, y^{k-\nu}
  \label{eq:brl} 
\shortintertext{and}
  \sum_{\nu=k}^{n} \binom{n+1}{\nu+1} \AP_{\nu-k}(x) \, y^{n-\nu}
  &= \sum_{\nu=k}^{n} \binom{\nu}{k} \AP_{n-\nu}(x+y) \, y^{\nu-k}.
  \label{eq:brr} 
\shortintertext{In particular,}
  \sum_{\nu=0}^{n} \binom{n+1}{\nu+1} \AP_\nu(x) \, y^{n-\nu}
  &= \sum_{\nu=0}^{n} \AP_\nu(x+y) \, y^{n-\nu}.
  \label{eq:brs} 
\end{align}
\end{theorem}

The binomial sum of the left-hand side of~\eqref{eq:brl} can be interpreted as an 
incomplete binomial sum compared to \eqref{eq:apt}, as well as shifted in the row, 
e.g., $\binom{k}{\nu} \mapsto \binom{n+1}{\nu}$.
Combining \eqref{eq:brl} and \eqref{eq:brr}, see Figure~\ref{fig:hs3}(b) below, 
gives a relationship between the last row $n+1$ and an upright V-shape, where the 
vertex with $\binom{n}{k}$ is counted twice.

\begin{figure}[H]
\begin{center}
\begin{subfigure}[b]{5cm}
\includegraphics[width=5cm]{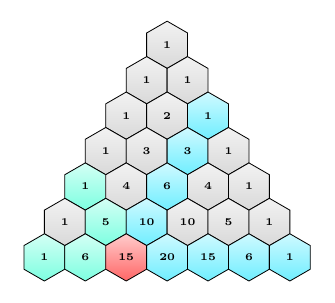}
\vspace*{-5ex}
\caption{}
\end{subfigure}
\begin{subfigure}[b]{5cm}
\includegraphics[width=5cm]{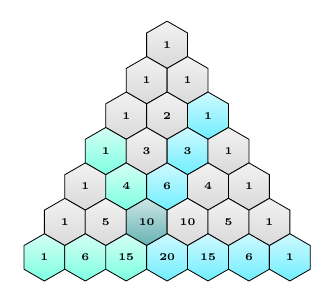}
\vspace*{-5ex}
\caption{}
\end{subfigure}
\caption{Asymmetric shapes of boomerangs in Pascal's triangle. 
(a) The red hexagonal area marks the omitted vertex.
The colored shapes on the left and right correspond to 
\eqref{eq:brl} with $(n,k)=(5,1)$ and \eqref{eq:brr} with $(n,k)=(5,2)$, respectively.
(b) Overlapping shapes of \eqref{eq:brl} and \eqref{eq:brr} with $(n,k)=(5,2)$, respectively.} 
\label{fig:hs3}
\end{center}
\end{figure}

By choosing $\AP_\nu(x) = x^\nu$ for $\nu \geq 0$, we obtain the following 
corollary at once.

\begin{corollary} \label{cor:main}
For $n \geq k \geq 0$, we have
\begin{align}
  \sum_{\nu=0}^{k} \binom{n+1}{\nu} x^{k-\nu} y^{\nu}
  &= \sum_{\nu=0}^{k} \binom{n-\nu}{k-\nu} (x+y)^{\nu} y^{k-\nu}
  \label{eq:brl2} 
\shortintertext{and}
  \sum_{\nu=k}^{n} \binom{n+1}{\nu+1} x^{\nu-k} y^{n-\nu}
  &= \sum_{\nu=k}^{n} \binom{\nu}{k} (x+y)^{n-\nu} y^{\nu-k}.
  \nonumber
\shortintertext{In particular,}
  \sum_{\nu=0}^{n} \binom{n+1}{\nu+1} x^{\nu} y^{n-\nu}
  &= \sum_{\nu=0}^{n} (x+y)^{\nu} y^{n-\nu}.
  \nonumber
\end{align}
\end{corollary}

Setting $y=1$ yields the following univariate identities.

\begin{corollary} \label{cor:main2}
For $n \geq k \geq 0$, we have
\begin{align}
  \sum_{\nu=0}^{k} \binom{n+1}{\nu} x^{k-\nu} 
  &= \sum_{\nu=0}^{k} \binom{n-\nu}{k-\nu} (x+1)^{\nu} 
  \label{eq:brl3} 
\shortintertext{and}
  \sum_{\nu=k}^{n} \binom{n+1}{\nu+1} x^{\nu-k} 
  &= \sum_{\nu=k}^{n} \binom{\nu}{k} (x+1)^{n-\nu}.
  \label{eq:brr3} 
\shortintertext{In particular,}
  \sum_{\nu=0}^{n} \binom{n+1}{\nu+1} x^{\nu}
  &= \sum_{\nu=0}^{n} (x+1)^{\nu}.
  \label{eq:brs3}
\end{align}
\end{corollary}

One step further, we finally arrive at the hockey-stick identities.

\begin{corollary} \label{cor:main3}
The hockey-stick identities follow from \eqref{eq:brl3} and \eqref{eq:brr3}
by taking the values
\begin{enumerate}
  \item $x = 0$ that implies \eqref{eq:hsl} and \eqref{eq:hsr},
  \item $x = -1$ that implies \eqref{eq:asl} and \eqref{eq:asr},
\end{enumerate}
respectively.
\end{corollary}

Finally, we obtain a kind of symmetry relations.

\begin{corollary} \label{cor:main4}
Let $n \geq k \geq 0$. We have a reflection relation by 
\begin{align}
  \sum_{\nu=0}^{k} \binom{n+1}{\nu} (-\tfrac{1}{2})^{k-\nu} 
  &= \sum_{\nu=0}^{k} \binom{n-\nu}{k-\nu} (\tfrac{1}{2})^\nu 
  \label{eq:brr4} 
\shortintertext{and an incomplete binomial sum by}
  \sum_{\nu=0}^{k} \binom{n+1}{\nu} 
  &= \sum_{\nu=0}^{k} \binom{n-\nu}{k-\nu} 2^\nu. 
  \label{eq:bri4}
\end{align}
In particular, involving the central binomial coefficient $\binom{2n}{n}$, we obtain
\begin{equation} \label{eq:brc4}
  \sum_{\nu=0}^{n} \binom{n+\nu}{\nu} 2^{n-\nu}
  = \sum_{\nu=0}^{n} \binom{n+\nu}{n} 2^{n-\nu} = 2^{2n}.
\end{equation}
\end{corollary}

\begin{figure}[H]
\begin{center}
\includegraphics[width=5cm]{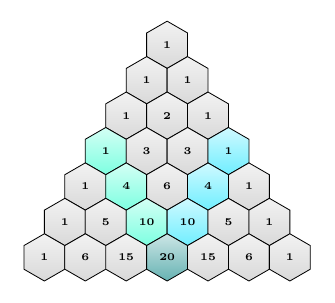}
\caption{Symmetric sums with the central binomial coefficient, which is counted
twice, as given by \eqref{eq:brc4}.} 
\label{fig:hs4}
\end{center}
\end{figure}

% ==============================================================================
\section{Applications}
\label{sec:appl}
% ==============================================================================

We will show several applications of Theorem~\ref{thm:main} and its corollaries,
see Figure~\ref{fig:diag}. The various identities that will be derived below are 
summarized in the appendix, see Tables~\ref{tab:term}~--~\ref{tab:term3}.

\tikzstyle{box} = [rectangle, rounded corners, text width=3cm,
minimum height=1cm, text centered, font=\small, draw=black, fill=cyan!10]

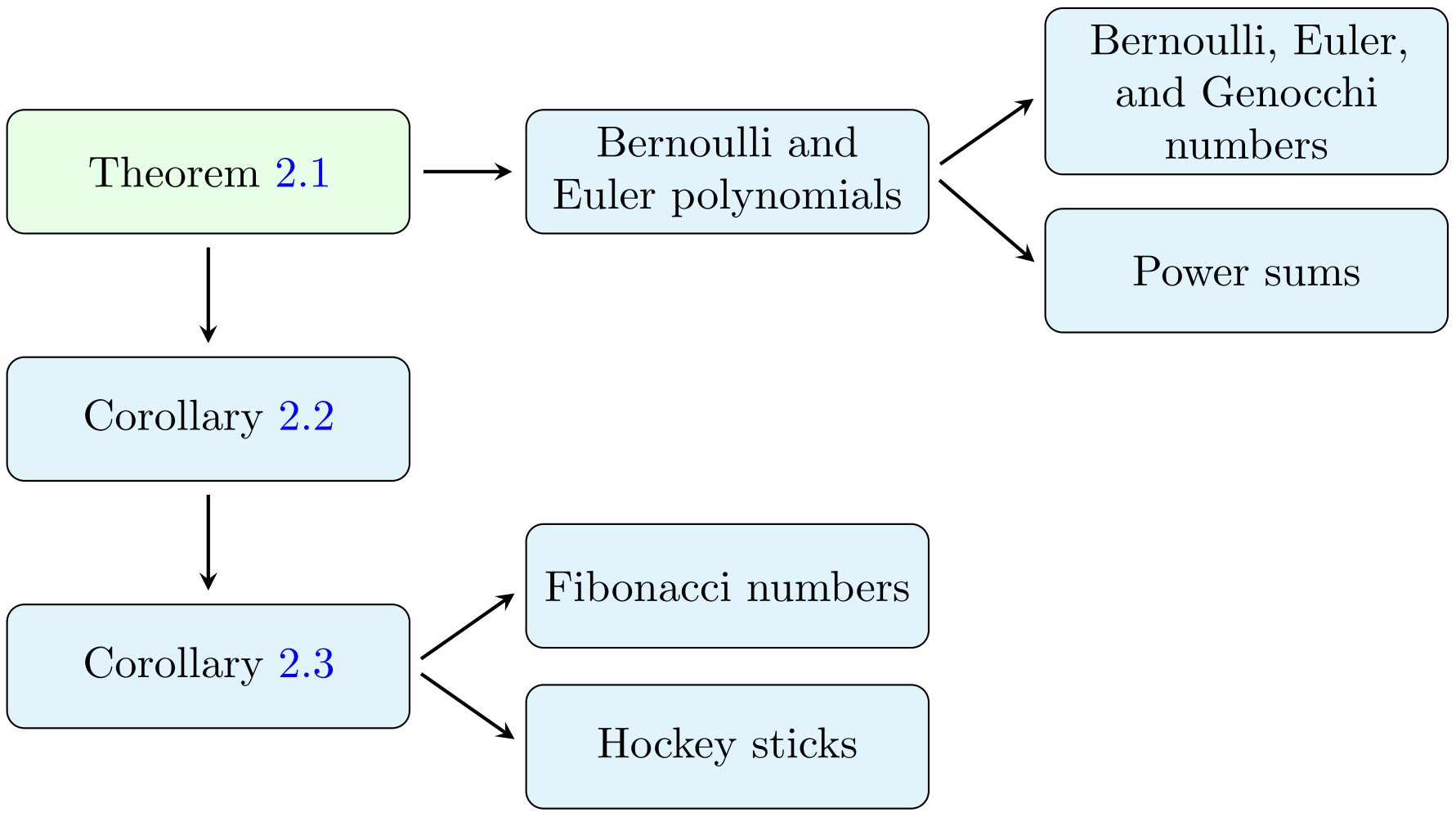
\begin{figure}[H]
\begin{center}
\begin{tikzpicture}[
    node/.style={font=\small},
    arrow/.style={->, >=stealth, black, thick, shorten >=3pt, shorten <=3pt}]

    \node[box, fill=green!10] (T1) at (0,4) {Theorem~\ref{thm:main}};
    \node[box] (C2) at (0,2) {Corollary~\ref{cor:main}};
    \node[box] (C3) at (0,0) {Corollary~\ref{cor:main2}};
    
    \node[box] (AP) at (4.2,4) {Bernoulli and Euler polynomials};
    \node[box] (AN) at (8.4,4.65) {Bernoulli, Euler, \\ and Genocchi\\ numbers};
    \node[box] (PS) at (8.4,3.20) {Power sums};
    
    \node[box] (F) at (4.2,0.65) {Fibonacci numbers};
    \node[box] (H) at (4.2,-0.65) {Hockey sticks};

    \draw[arrow] (T1) -- (AP);
    \draw[arrow] (T1) -- (C2);
    \draw[arrow] (AP.east) -- (AN.west);
    \draw[arrow] (AP.east) -- (PS.west);
    \draw[arrow] (C2) -- (C3);
    \draw[arrow] (C3.east) -- (H.west);
    \draw[arrow] (C3.east) -- (F.west);
\end{tikzpicture}
\medskip
\caption{Applications}
\label{fig:diag}
\end{center}
\end{figure}

The obtained identities mainly contain those binomial coefficients
that arise from the boomerang identities~\eqref{eq:brl} and \eqref{eq:brs}. 
We shall use increasing indices and powers in formulas, e.g., $\AP_\nu(x)$ and $x^\nu$,
respectively. Correspondingly, we have for $n \geq k \geq 0$ the identities
\begin{align}
  \sum_{\nu=0}^{k} \binom{n+1}{k-\nu} \, \RT_\nu
  &= \sum_{\nu=0}^{k} \binom{n-\nu}{k-\nu} \, \CT_\nu \;+\; \ST_k 
  \label{eq:btl} 
\shortintertext{and the limiting cases}  
  \sum_{\nu=0}^{n} \binom{n+1}{\nu+1} \, \RT_\nu
  &= \sum_{\nu=0}^{n} \CT_\nu \;+\; \ST_n, 
  \label{eq:bts} 
\end{align}
where $\RT$, $\CT$, and $\ST$ represent terms in the rows and
columns of Pascal's triangle, as well as, supplementary terms, respectively.
For simplicity, we omit the equivalent cases that would follow from 
\eqref{eq:brr} by the symmetry, see Figure~\ref{fig:hs5} below.

\begin{figure}[H]
\begin{center}
\begin{subfigure}[b]{5cm}
\includegraphics[width=5cm]{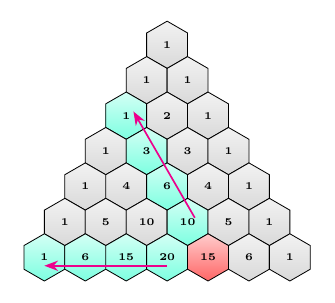}
\vspace*{-5ex}
\caption{}
\end{subfigure}
\begin{subfigure}[b]{5cm}
\includegraphics[width=5cm]{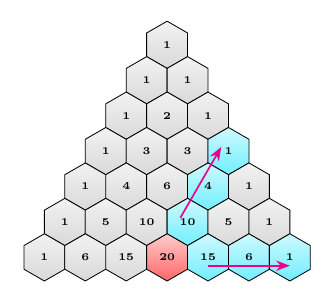}
\vspace*{-5ex}
\caption{}
\end{subfigure}
\begin{subfigure}[b]{5cm}
\includegraphics[width=5cm]{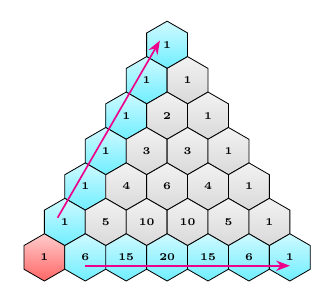}
\vspace*{-5ex}
\caption{}
\end{subfigure}
\caption{Colored shapes of boomerangs in Pascal's triangle, where (a) and (c) 
correspond to \eqref{eq:btl} and \eqref{eq:bts}, respectively. Case (b) is omitted.
The red hexagonal areas mark the omitted vertices, from which the arrows point 
in the direction of the increasing indices of the terms.} 
\label{fig:hs5}
\end{center}
\end{figure}

As usual, empty sums are defined as $0$.
Let denote ${\delta_c = 1}$ if condition $c$ is true, $0$~otherwise.
(Kronecker delta meets Iverson bracket.)

% ------------------------------------------------------------------------------
\subsection{Fibonacci identities}
% ------------------------------------------------------------------------------

The sequence of the Fibonacci numbers for $n \geq 0$ begins
\[
  \FN_n = 0, 1, 1, 2, 3, 5, 8, 13, 21, 34, 55, 89, 144, 233, \ldots,
\]
obeying the recurrence ${\FN_n = \FN_{n-1} + \FN_{n-2}}$ for ${n \geq 2}$.
See sequence \seqnum{A000045} in the OEIS~\cite{OEIS}.
Extending to negative indices, it follows that 
${\FN_{-n} = (-1)^{n+1} \FN_n}$ for ${n \geq 0}$.

\begin{theorem} \label{thm:fib}
For $n \geq k \geq 0$, we have
\begin{align}
  \sum_{\nu=0}^{k} \binom{n+1}{k-\nu} \FN_\nu 
  &= \sum_{\nu=0}^{k} \binom{n-\nu}{k-\nu} \FN_{2\nu} 
  \label{eq:tf1} 
\shortintertext{and}
  \sum_{\nu=0}^{k} \binom{n+1}{k-\nu} (-1)^\nu \FN_{2\nu} 
  &= \sum_{\nu=0}^{k} \binom{n-\nu}{k-\nu} (-1)^\nu \FN_\nu. 
  \label{eq:tf2}
\end{align}
Moreover, there are the recurrence relations
\begin{equation} \label{eq:tf3}
  \sum_{\nu=0}^{k} \left( (-1)^\nu \binom{n+1}{k-\nu} + \binom{n-\nu}{k-\nu} \right) \FN_\nu = 0 
\end{equation}
and with even indices that
\begin{equation} \label{eq:tf4}
  \sum_{\nu=0}^{k} \left( (-1)^\nu \binom{n+1}{k-\nu} + \binom{n-\nu}{k-\nu} \right) 
  \left( \FN_{2\nu} - 2 \delta_{2 \mid \nu} \, \FN_\nu \right) = 0.
\end{equation}
\end{theorem}

For $n=k$, we obtain the following well-known identities (see sequences 
\seqnum{A027941} and \seqnum{A119282}, respectively).

\begin{corollary} \label{cor:fib}
For $n \geq 0$, we have
\[
  \FN_{2n+1} - 1 = \sum_{\nu=0}^{n} \binom{n+1}{\nu+1} \FN_\nu
  = \sum_{\nu=0}^{n} \FN_{2\nu},
\]
as well as,
\[
  (-1)^n \FN_{n-1} - 1 = \sum_{\nu=0}^{n} \binom{n+1}{\nu+1} (-1)^\nu \FN_{2\nu} 
  = \sum_{\nu=0}^{n} (-1)^\nu \FN_\nu.
\]
\end{corollary}

% ------------------------------------------------------------------------------
\subsection{Appell and Bernoulli identities}
% ------------------------------------------------------------------------------

In this subsection, let $(\AP_\nu(x))_{\nu \geq 0}$ be any Appell sequence that
satisfies the reflection relation~\eqref{eq:apr}.

\begin{theorem} \label{thm:app}
For $n \geq k \geq 0$, we have
\begin{align}
  \sum_{\nu=0}^{k} \binom{n+1}{k-\nu} \AP_\nu(x) 
  &= \sum_{\nu=0}^{k} \binom{n-\nu}{k-\nu} (-1)^\nu \AP_\nu(-x). 
  \label{eq:app-ref1} 
\shortintertext{In particular,}
  \sum_{\nu=0}^{n} \binom{n+1}{\nu+1} \AP_\nu(x) 
  &= \sum_{\nu=0}^{n} (-1)^\nu \AP_\nu(-x).
  \label{eq:app-ref2}
\end{align}
Taking $x=0$ and applying the Bernoulli and Euler polynomials $\BN_\nu(x)$ and 
$\EN_\nu(x)$, we can replace $\AP_\nu(0)$ by the values
\begin{enumerate}
\item $\BN_\nu(0) = \BN_\nu$,
\item $\EN_\nu(0) \hspace*{0.2ex} = \GD_{\nu+1}$,
\end{enumerate}
where ${\BN_\nu \in \QQ}$ are the Bernoulli numbers, 
${\GN_\nu = 2(1-2^\nu) \BN_\nu \in \ZZ}$ are the Genocchi numbers, 
and $\GD_{\nu+1} = \GN_{\nu+1}/(\nu+1) \in \QQ$ are the divided Genocchi numbers, 
respectively.
\end{theorem}

For the Bernoulli polynomials and numbers, this yields the following results.

\begin{corollary} \label{cor:app-bp}
For $n \geq k \geq 0$, we have
\begin{align}
  \sum_{\nu=0}^{k} \left( \binom{n+1}{k-\nu} - \binom{n-\nu}{k-\nu} \right) \BN_\nu(x) 
  &= \sum_{\nu=1}^{k} \binom{n-\nu}{k-\nu} \nu x^{\nu-1} 
  \label{eq:app-bp1} 
\shortintertext{and}
  \sum_{\nu=0}^{k} \left( \binom{n+1}{k-\nu} - \binom{n-\nu}{k-\nu} \right) \BN_\nu 
  &= \binom{n-1}{k-1}. 
  \label{eq:app-bp2} 
\shortintertext{In particular,}
  \sum_{\nu=0}^{n} \left( \binom{n+1}{\nu+1} - 1 \right) \BN_\nu(x) 
  &= \sum_{\nu=1}^{n} \nu x^{\nu-1} 
  \label{eq:app-bp3} 
\shortintertext{and}
  \sum_{\nu=0}^{n} \left( \binom{n+1}{\nu+1} - 1 \right) \BN_\nu
  &= \delta_{n > 0}. 
  \label{eq:app-bp4} 
\end{align}
\end{corollary}

This provides an infinite family of recurrences for the Bernoulli numbers.

\begin{corollary} \label{cor:app-bn}
For $n \geq k \geq 3$ with odd $k$, we have the recurrences
\begin{align}
  \BN_{k-1} &= \frac{1}{2k} \left( \binom{n-1}{k-3} 
  - 2 \sum_{\substack{\nu=2\\2 \mids \nu}}^{k-3} 
  \left( \binom{n+1}{k-\nu} - \binom{n-\nu}{k-\nu} \right) \BN_\nu \right).
  \label{eq:app-bn1} 
\shortintertext{In particular, we have for even $n \geq 2$ that}
  \BN_n &= \frac{1}{2(n+1)} \left( \binom{n}{2} 
  - 2 \sum_{\substack{\nu=2\\2 \mids \nu}}^{n-2} 
  \left( \binom{n+2}{\nu+1} - 1 \right) \BN_\nu \right).
  \label{eq:app-bn2} 
\end{align}
\end{corollary}

% ------------------------------------------------------------------------------
\subsection{Euler polynomials and Genocchi numbers}
% ------------------------------------------------------------------------------

From Theorem~\ref{thm:app}, we deduce similar results for the Euler polynomials 
and divided Genocchi numbers.

\begin{corollary} \label{cor:app-ep}
For $n \geq k \geq 0$, we have
\begin{align}
  \sum_{\nu=0}^{k} \left( \binom{n+1}{k-\nu} + \binom{n-\nu}{k-\nu} \right) \EN_\nu(x)
  &= 2 \sum_{\nu=0}^{k} \binom{n-\nu}{k-\nu} x^\nu 
  \label{eq:app-ep1} 
\shortintertext{and}
  \sum_{\nu=0}^{k} \left( \binom{n+1}{k-\nu} + \binom{n-\nu}{k-\nu} \right) \GD_{\nu+1} 
  &= 2 \binom{n}{k}. 
  \label{eq:app-ep2} 
\shortintertext{In particular,}
  \sum_{\nu=0}^{n} \left( \binom{n+1}{\nu+1} + 1 \right) \EN_\nu(x) 
  &= 2 \sum_{\nu=0}^{n} x^\nu 
  \label{eq:app-ep3} 
\shortintertext{and}
  \sum_{\nu=0}^{n} \left( \binom{n+1}{\nu+1} + 1 \right) \GD_{\nu+1}
  &= 2. \label{eq:app-ep4} 
\end{align}
\end{corollary}

Correspondingly, we obtain an infinite number of recurrence relations for the 
divided Genocchi numbers.

\begin{corollary} \label{cor:app-gn}
For $n \geq k \geq 1$ and odd $k$, we have the recurrences
\begin{equation} \label{eq:app-gn1}
  \GD_{k+1} = - \frac{1}{2} \left( \binom{n}{k-1} 
  + \sum_{\substack{\nu=1\\2 \nmids \nu}}^{k-2} 
  \left( \binom{n+1}{k-\nu} + \binom{n-\nu}{k-\nu} \right) \GD_{\nu+1} \right).
\end{equation}
In particular, we have for even $n \geq 2$ that
\begin{equation} \label{eq:app-gn2}
  \GD_n = - \frac{1}{2} \left( \binom{n-1}{1} 
  + \sum_{\substack{\nu=2\\2 \mids \nu}}^{n-2} 
  \left( \binom{n}{\nu} + 1 \right) \GD_\nu \right).
\end{equation}
\end{corollary}

% ------------------------------------------------------------------------------
\subsection{Genocchi numbers again}
% ------------------------------------------------------------------------------

Theorem~\ref{thm:app} also implies the following identities of the Genocchi numbers.

\begin{corollary} \label{cor:app-dgn}
For $n \geq k \geq 0$, we have
\begin{align}
  \sum_{\nu=0}^{k} \binom{n+1}{k-\nu} \GN_\nu
  &= \sum_{\nu=0}^{k} \binom{n-\nu}{k-\nu} (-1)^{\nu+1} \GN_\nu. 
  \label{eq:app-dgn1}
\shortintertext{In particular,}
  \sum_{\nu=0}^{n} \binom{n+1}{\nu+1} \GN_\nu 
  &= \sum_{\nu=0}^{n} (-1)^{\nu+1} \GN_\nu.
  \label{eq:app-dgn2}
\end{align}
\end{corollary}

\begin{corollary} \label{cor:app-dgn2}
For $n \geq k \geq 0$, we have
\begin{align}
  \sum_{\nu=0}^{k} \left( \binom{n+1}{k-\nu} + \binom{n-\nu}{k-\nu} \right) \GN_\nu 
  &= 2 \binom{n-1}{k-1}. 
  \label{eq:app-dgn3}
\shortintertext{In particular,}
  \sum_{\nu=0}^{n} \left( \binom{n+1}{\nu+1} + 1 \right) \GN_\nu
  &= 2\delta_{n > 0}.
  \label{eq:app-dgn4}
\end{align}
\end{corollary}

This also provides an infinite number of recurrence relations for the Genocchi numbers.

\begin{corollary} \label{cor:app-dgn3} 
For $n \geq k \geq 2$ and even $k$, we have
\begin{equation} \label{eq:app-dgn5}
  \GN_k = - \frac{1}{2} \left( \binom{n+1}{k-1} - \binom{n-1}{k-1} 
  + \sum_{\substack{\nu=2\\2 \nmids \nu}}^{k-2} 
  \left( \binom{n+1}{k-\nu} + \binom{n-\nu}{k-\nu} \right) \GN_\nu \right).
\end{equation}
In particular, we have for even $n \geq 2$ that
\begin{equation} \label{eq:app-dgn6}
  \GN_n = - \frac{1}{2} \left( \binom{n+1}{2} - 1 
  + \sum_{\substack{\nu=2\\2 \mids \nu}}^{n-2} 
  \left( \binom{n+1}{\nu+1} + 1 \right) \GN_\nu \right).
\end{equation}
\end{corollary}

% ------------------------------------------------------------------------------
\subsection{Mixed identities}
% ------------------------------------------------------------------------------

The reciprocal Appell polynomials are defined by
\[
  \AR_n(x) = x^n \AP_n(x^{-1}),
\] 
where the coefficients of $\AR_n(x)$ are listed in reverse order and therefore 
defined for all $x \in \CC$.

\begin{theorem} \label{thm:mix}
For $n \geq k \geq 0$, we have
\begin{align}
  \sum_{\nu=0}^{k} \binom{n+1}{k-\nu} x^\nu \AP_\nu(0)
  &= \sum_{\nu=0}^{k} \binom{n-\nu}{k-\nu} \AR_\nu(x)
  \label{eq:mix-ap1} 
\shortintertext{and}
  \sum_{\nu=0}^{k} \binom{n+1}{k-\nu} (-1)^\nu \AR_\nu(x)
  &= \sum_{\nu=0}^{k} \binom{n-\nu}{k-\nu} (-x)^\nu \AP_\nu(0). 
  \label{eq:mix-ap2} 
\shortintertext{In particular,}
  \sum_{\nu=0}^{n} \binom{n+1}{\nu+1} x^\nu \AP_\nu(0) 
  &= \sum_{\nu=0}^{n} \AR_\nu(x) 
  \label{eq:mix-ap3} 
\shortintertext{and}
  \sum_{\nu=0}^{n} \binom{n+1}{\nu+1} (-1)^\nu \AR_\nu(x) 
  &= \sum_{\nu=0}^{n} (-x)^\nu \AP_\nu(0).
  \label{eq:mix-ap4} 
\end{align}
\end{theorem}

\begin{corollary} \label{cor:mix-hs}
Assuming $\AP_0(0) \neq 0$ and setting $x=0$ in \eqref{eq:mix-ap1} and \eqref{eq:mix-ap2} 
imply both kinds of the hockey-stick identities \eqref{eq:hsl} and \eqref{eq:asl}, 
respectively. 
\end{corollary}

Now, let $(\AP_\nu(x))_{\nu \geq 0}$ be any Appell sequence that satisfies 
the reflection relation~\eqref{eq:apr}. Special relations are given for $x=2$.

\begin{corollary} \label{cor:mix}
For $n \geq k \geq 0$, we have
\begin{align}
  \sum_{\nu=0}^{k} \binom{n+1}{k-\nu} 2^\nu \AP_\nu(0)
  &= \sum_{\substack{\nu=0\\2 \mids \nu}}^{k} \binom{n-\nu}{k-\nu} \AR_\nu(2)
  \label{eq:mix2-ap1} 
\shortintertext{and}
  \sum_{\substack{\nu=0\\2 \mids \nu}}^{k} \binom{n+1}{k-\nu} \AR_\nu(2)
  &= \sum_{\nu=0}^{k} \binom{n-\nu}{k-\nu} (-2)^\nu \AP_\nu(0) 
  \label{eq:mix2-ap2} \\
  &= \sum_{\nu=0}^{k} \binom{n-\nu}{k-\nu} (-1)^\alpha \, 2^\nu \AP_\nu(0) \;+\; \ST_k. 
  \label{eq:mix2-ap3} 
\shortintertext{In particular,}
  \sum_{\nu=0}^{n} \binom{n+1}{\nu+1} 2^\nu \AP_\nu(0) 
  &= \sum_{\substack{\nu=0\\2 \mids \nu}}^{n} \AR_\nu(2)
  \label{eq:mix2-ap4} 
\shortintertext{and}
  \sum_{\substack{\nu=0\\2 \mids \nu}}^{n} \binom{n+1}{\nu+1} \AR_\nu(2) 
  &= \sum_{\nu=0}^{n} (-2)^\nu \AP_\nu(0) 
  \label{eq:mix2-ap5} \\
  &= \sum_{\nu=0}^{n} (-1)^\alpha \, 2^\nu \AP_\nu(0) \;+\; \ST_n.
  \label{eq:mix2-ap6} 
\end{align}
Applying the Bernoulli and Euler polynomials $\BN_\nu(x)$ and $\EN_\nu(x)$,
the values can be replaced by
\begin{enumerate}
\item $\BN_\nu(0) = \BN_\nu$, $\BR_\nu(2) = \BN_\nu + \tfrac{1}{2} \GN_\nu$, 
  $\alpha = 0$, $\ST_k = 2\binom{n-1}{k-1}$, and $\ST_n = 2\delta_{n > 0}$,
\item $\EN_\nu(0) = \GD_{\nu+1}$, $\ER_\nu(2) = \EN_\nu$, 
  $\alpha = 1$, $\ST_k = 2\binom{n}{k}$, and $\ST_n = 2$,
\end{enumerate}
where ${\EN_\nu \in \ZZ}$ are the Euler numbers, respectively.
\end{corollary}

\begin{corollary} \label{cor:mix-en}
For $n \geq 0$, we have
\begin{align}
  \sum_{\nu=0}^{n} \EN_\nu 
  &= \frac{1}{2} \sum_{\nu=1}^{n+1} \binom{n+1}{\nu} 2^\nu \GD_{\nu}
  \label{eq:en-sum} 
\shortintertext{and}
  \EN_n &= \frac{1}{2(n+1)} \sum_{\nu=1}^{n+1} \binom{n+1}{\nu} 2^\nu \GN_{\nu}.
  \label{eq:en-sum2}
\end{align}
\end{corollary}

\begin{example}
Using Corollaries~\ref{cor:mix} and \ref{cor:mix-en}, we evaluate by Table~\ref{tab:ap-val} 
two different sum identities of the Euler and divided Genocchi numbers. 
First, we obtain for $n=9$ and $k=7$ that
\begin{align*}
  \sum_{\nu=0}^{k} \binom{n+1}{k-\nu} \EN_\nu
  &= \binom{10}{7} - \binom{10}{5} + \binom{10}{3} 5 - \binom{10}{1} 61 = -142 
\shortintertext{and}
  \sum_{\nu=0}^{k} \binom{n-\nu}{k-\nu} (-2)^\nu \GD_{\nu+1} 
  &= \binom{9}{7} (-2)^0 - \binom{8}{6} (-2)^1 \tfrac{1}{2}
  + \binom{6}{4} (-2)^3 \tfrac{1}{4} \\
  &\hspace*{2.5ex} - \binom{4}{2} (-2)^5 \tfrac{1}{2}
  + \binom{2}{0} (-2)^7 \tfrac{17}{8} = -142. 
\end{align*}
Second, we have for $n=8$ that
\begin{align*}
  \sum_{\nu=0}^{n} \EN_\nu &= 1 - 1 + 5 - 61 + 1385 = 1329 
\shortintertext{and}
  \frac{1}{2} \sum_{\nu=1}^{n+1} \binom{n+1}{\nu} 2^\nu \GD_{\nu}
  &= \binom{9}{1} 2^0 - \binom{9}{2} 2^1 \tfrac{1}{2}
  + \binom{9}{4} 2^3 \tfrac{1}{4} \\
  &\hspace*{2.5ex} - \binom{9}{6} 2^5 \tfrac{1}{2}
  + \binom{9}{8} 2^7 \tfrac{17}{8} = 1329.
\end{align*}
\end{example}

% ------------------------------------------------------------------------------
\subsection{Power sum identities}
% ------------------------------------------------------------------------------

Define the power sums
\[
  S_n(m) = \sum_{j=0}^{m-1} j^n 
\]
for $m, n \geq 0$, which are given by polynomials $S_n(x)$ of degree $n+1$.
Further define the differences 
\[
  S_n(x_1,x_2) = S_n(x_2) - S_n(x_1).
\]

\begin{theorem} \label{thm:sum}
For $n \geq k \geq 0$, we have
\begin{align*}
  \sum_{\nu=0}^{k} \binom{n+1}{k-\nu} S_\nu(x_1,x_2) \, y^{k-\nu} 
  &= \sum_{\nu=0}^{k} \binom{n-\nu}{k-\nu} S_\nu(x_1+y,x_2+y) \, y^{k-\nu}. 
\shortintertext{In particular,}
  \sum_{\nu=0}^{n} \binom{n+1}{\nu+1} S_\nu(x_1,x_2) \, y^{k-\nu} 
  &= \sum_{\nu=0}^{n} S_\nu(x_1+y,x_2+y) \, y^{k-\nu}.
\end{align*}
\end{theorem}

\begin{corollary} \label{cor:sum}
For $n \geq k \geq 0$, we have
\begin{align*}
  \sum_{\nu=0}^{k} \left( \binom{n+1}{k-\nu} - \binom{n-\nu}{k-\nu} \right) S_\nu(x) 
  &= \sum_{\nu=1}^{k} \binom{n-\nu}{k-\nu} x^\nu. 
\shortintertext{In particular,}
  \sum_{\nu=0}^{n} \left( \binom{n+1}{\nu+1} - 1 \right) S_\nu(x) 
  &= \sum_{\nu=1}^{n} x^\nu. 
\end{align*}
\end{corollary}

% ==============================================================================
\section{Preliminaries}
\label{sec:prel}
% ==============================================================================

The Fibonacci numbers can be defined for $n \geq 0$ by
\begin{equation} \label{eq:fib}
  \FN_n = \frac{\phi^n - (1-\phi)^n}{\sqrt{5}}
\end{equation}
with the golden ratio $\phi = (1 + \sqrt{5})/2$. Furthermore, we have the 
well-known relations $\phi^2 = \phi + 1$ and $(1-\phi)^2 = 2 - \phi$
(cf.~\cite[Sec.\,6.6, p.\,299]{GKP:1994}).
For polynomials $f \in \CC[x]$, we shall evaluate 
\begin{equation} \label{eq:f-phi}
  f(\FN) = \frac{f(\phi) - f(1-\phi)}{\sqrt{5}}
\end{equation}
to obtain sums connected with the Fibonacci numbers, where $f(\FN)$ is to be 
understood in \emph{umbral sense} such that ${\FN^\nu \mapsto \FN_\nu}$.
In this context, we consider the following mappings that exchange shifted 
monomials with the Fibonacci numbers.

\begin{lemma} \label{lem:map}
Regarding \eqref{eq:f-phi}, we have for shifted monomials for $\nu \geq 0$
the mappings 
\[
  x^\nu \mapsto \FN_\nu, \quad (x+1)^\nu \mapsto \FN_{2\nu}, 
  \quad (x-1)^\nu \mapsto (-1)^{\nu+1} \FN_\nu,
  \quad (1-x)^\nu \mapsto -\FN_\nu.
\]
\end{lemma}

\begin{proof}
For each monomial $x^\nu$ and shifted ones, we have to compute 
$(x^\nu - (1-x)^\nu)/\sqrt{5} \valueat{x=\phi}$ in view of \eqref{eq:f-phi}
and to rearrange terms to coincide with \eqref{eq:fib}. This leads to the 
following table with computed terms according to $\phi$:
\begin{center} 
\begin{tabular}{cllr@{$\;\mapsto\;$}l}
\toprule
  $x$   & $\phi$ & $1-\phi$ & $x^\nu$ & $\FN_\nu$ \\
  $x+1$ & $\phi+1 = \phi^2$ & $2-\phi = (1-\phi)^2$ & $(x+1)^\nu$ & $\FN_{2\nu}$ \\
  $x-1$ & $\phi-1 = -(1-\phi)$ & $-\phi$ & $(x-1)^\nu$ & $(-1)^{\nu+1} \FN_\nu$ \\
  $1-x$ & $1-\phi$ & $\phi$ & $(1-x)^\nu$ & $-\FN_\nu$. \\
\bottomrule
\end{tabular}
\end{center} 
This implies the results.
\end{proof}

\begin{proposition} \label{prop:ident}
Let $f,g \in \CC[x]$ of degree $n \geq 0$, where
\[
  f(x) = \sum_{\nu=0}^{n} a_\nu \, x^\nu \andq 
  g(x) = \sum_{\nu=0}^{n} b_\nu \, x^\nu.
\]
If $f(x) = g(x+1)$, then we have the identities
\begin{align}
  \sum_{\nu=0}^{n} a_\nu \, \FN_\nu 
  &= \sum_{\nu=0}^{n} b_\nu \, \FN_{2\nu} 
  \label{eq:fib-id1} 
\shortintertext{and}
  \sum_{\nu=0}^{n} (-1)^\nu a_\nu \, \FN_{2\nu} 
  &= \sum_{\nu=0}^{n} (-1)^\nu b_\nu \, \FN_\nu.
  \label{eq:fib-id2}
\end{align}
Moreover, there are the recurrence relations
\begin{align}
  \sum_{\nu=0}^{n} \left( (-1)^{\nu} a_\nu + b_\nu \right) \FN_\nu &= 0
  \label{eq:fib-id3} \\
\shortintertext{and}
  \sum_{\nu=0}^{n} \left( (-1)^{\nu} a_\nu + b_\nu \right) 
  (\FN_{2\nu} - (-1)^\nu \FN_\nu) &= 0.
  \label{eq:fib-id4} 
\end{align}
\end{proposition}

\begin{proof}
We apply Lemma~\ref{lem:map} for $f$ and $g$ using the pairs ${(x,x+1)}$, 
${(-(x+1),-x)}$, and ${(x-1,x)}$, which provide \eqref{eq:fib-id1}, 
\eqref{eq:fib-id2}, and \eqref{eq:fib-id3}, respectively. 
Subtracting \eqref{eq:fib-id1} from \eqref{eq:fib-id2} finally yields 
\eqref{eq:fib-id4}, finishing the proof.
\end{proof}

The Bernoulli and Euler polynomials, as well as, the Euler and Genocchi numbers 
(cf.~\cite[Sec.\,1.14, pp.\,48--49]{Comtet:1974}) are defined by the exponential 
generating functions
\begin{alignat*}{5}
  \frac{t e^{xt}}{e^t - 1} &= \sum_{n \geq 0} \BN_n(x) \frac{t^n}{n!} \quad 
  && (|t| < 2\pi), \quad &
  \frac{2e^{xt}}{e^t+1} &= \sum_{n \geq 0} \EN_n(x) \frac{t^n}{n!} \quad 
  && (|t| < \pi), \\
  \frac{2 e^t}{e^{2t}+1} &= \sum_{n \geq 0} \EN_n \frac{t^n}{n!} \quad 
  && (|t| < \tfrac{\pi}{2}), \quad & 
  \frac{2t}{e^t+1} &= \sum_{n \geq 0} \GN_n \frac{t^n}{n!} \quad 
  && (|t| < \pi),
\end{alignat*}
respectively. The polynomials $\BN_n(x)$ and $\EN_n(x)$ are monic polynomials 
of degree $n$ for ${n \geq 0}$. There are the known relations that
\begin{equation} \label{eq:ap-val}
  \BN_n(0) = \BN_n \in \QQ, \quad
  \EN_n(0) = \GD_{n+1} = \tfrac{\GN_{n+1}}{n+1} \in \QQ, \andq
  \GN_n = 2(1-2^n) \BN_n \in \ZZ.
\end{equation}

The Bernoulli numbers obey the recurrence relation
\begin{equation} \label{eq:bn-rec}
  \sum_{\nu=0}^{n} \binom{n+1}{\nu} \BN_\nu = \delta_{n = 0},
\end{equation}
which also follows from \eqref{eq:apt} and \eqref{eq:apr} for $\BN_n(x)$ and 
setting $(x,y)=(0,1)$. The first few values of the Bernoulli, Euler, 
and (divided) Genocchi numbers are given in Table~\ref{tab:ap-val}.
For odd $n \geq 3$, the values of $\BN_n$, $\EN_n$, $\GN_n$, and $\GD_n$ vanish.
Moreover, there are the symmetry relations that
\begin{align}
  (-1)^n \, \BN_n(-x) &= \BN_n(x) + 
  \begin{cases}
    0, & \text{if $n=0$}, \\
    n x^{n-1}, & \text{otherwise}, 
  \end{cases}
  \label{eq:bp-sym} \\
  (-1)^n \, \EN_n(-x) &= - \EN_n(x) + 2 x^n. \label{eq:ep-sym}
\end{align}

\begin{table}[H] \small
\setstretch{1.25}
\begin{center}
\begin{tabular}{c|*{6}{r}}
  \toprule
  $n$ & $0$ & $1$ & $2$ & $4$ & $6$ & $8$ \\\hline
  $\BN_n$ & $1$ & $-\tfrac{1}{2}$ & $\tfrac{1}{6}$ & 
  $-\tfrac{1}{30}$ & $\tfrac{1}{42}$ & $-\tfrac{1}{30}$ \\
  $\EN_n$ & $1$ & $0$ & $-1$ & $5$ & $-61$ & $1385$ \\
  $\GN_n$ & $0$ & $1$ & $-1$ & $1$ & $-3$ & $17$ \\
  $\GD_n$ & & $1$ & $-\tfrac{1}{2}$ & $\tfrac{1}{4}$ & $-\tfrac{1}{2}$ & $\tfrac{17}{8}$ \\
  \bottomrule
\end{tabular}

\caption{First few values of $\BN_n$, $\EN_n$, $\GN_n$, and $\GD_n$.}
\label{tab:ap-val}
\end{center}
\end{table}

The reciprocal Appell polynomials are defined by
\begin{equation} \label{eq:ap-rec-def}
  \AR_n(x) = x^n \AP_n(x^{-1}).
\end{equation} 
Since the coefficients of $\AR_n(x)$ are listed in reverse order, we have 
\begin{equation} \label{eq:ap-rec0}
  \AR_n(0) = \AP_0(0).
\end{equation}
If $(\AP_\nu(x))_{\nu \geq 0}$ satisfies the reflection relation~\eqref{eq:apr}, then
\begin{equation} \label{eq:ap-rec1}
  \AR_n(2) = 0 \quad (n \geq 1 \text{ odd}).
\end{equation}
In particular, we have the special values
\begin{equation} \label{eq:ap-rec2}
  \BR_n(2) = \BN_n + \tfrac{1}{2} \GN_n \andq 
  \ER_n(2) = \EN_n.
\end{equation}
These properties above are discussed by the author \cite{Kellner:2021}.
The power sums
\[
  S_n(m) = \sum_{j=0}^{m-1} j^n 
\]
for $m, n \geq 0$ are given by
\begin{equation} \label{eq:sn-int}
  S_n(x) = \int_{0}^{x} \! \BN_n(t) \, dt = \frac{\BN_{n+1}(x) - \BN_{n+1}}{n+1}.
\end{equation}
Define the linear forward difference operator as $\Delta f(x) = f(x+1)-f(x)$,
so rule \eqref{eq:bin-add} becomes $\Delta \binom{n}{k+1} = \binom{n}{k}$.
In particular, \eqref{eq:bp-sym} is equivalent to
\[
  \Delta \BN_{n}(x) =
  \begin{cases}
    0, & \text{if $n=0$}, \\
    n x^{n-1}, & \text{otherwise}. 
  \end{cases}
\]
Consequently,
\begin{equation} \label{eq:sn-delta}
  \Delta S_n(x) = x^n.
\end{equation}

The well-known Vandermonde convolution for $x, y \in \CC$ and $n \geq 0$ reads
\begin{align}
  \sum_{\nu=0}^{n} \binom{x}{\nu} \binom{y}{n-\nu} &= \binom{x+y}{n}. 
  \nonumber 
\shortintertext{The corresponding \emph{upside-down} convolution 
(cf.~\cite[Sec.\,5.1, Eq.\,(5.26), p.\,169]{GKP:1994}) 
for integers $x, y \geq 0$ then reads}
  \sum_{\nu=0}^{n} \binom{\nu}{x} \binom{n-\nu}{y} &= \binom{n+1}{x+y+1}.
  \label{eq:ud-conv}
\end{align}
This identity can also be found in Gould's tables of combinatorial identities of 
1972 as \cite[Tab.\,3, Eq.~(3.3), p.~22]{Gould:1972}

% ==============================================================================
\section{Proofs}
\label{sec:proofs}
% ==============================================================================

% ------------------------------------------------------------------------------
\subsection{Main theorem}
% ------------------------------------------------------------------------------

\begin{proofcref}{thm:main}
Let $(\AP_\nu(x))_{\nu \geq 0}$ be an Appell sequence.
Let $n \geq k \geq 0$. We have to show that \eqref{eq:brl} holds, namely
\begin{equation} \label{eq:ap-loc}
  \sum_{\nu=0}^{k} \binom{n+1}{\nu} \AP_{k-\nu}(x) \, y^{\nu}
  = \sum_{\nu=0}^{k} \binom{n-\nu}{k-\nu} \AP_\nu(x+y) \, y^{k-\nu}.
\end{equation}
Applying \eqref{eq:bin-sym} and \eqref{eq:apt}, we can rewrite the right-hand side 
of \eqref{eq:ap-loc}, denoted by $R$, as
\[
  R = \sum_{\nu=0}^{k} \binom{n-\nu}{n-k} \sum_{\ell=0}^{\nu} 
  \binom{\nu}{\ell} \AP_\ell(x) \, y^{k-\ell}.
\]
We can increase the limits of the outer and inner sum, because then the corresponding 
binomial coefficients $\binom{n-\nu}{n-k}$ and $\binom{\nu}{\ell}$ vanish, respectively. 
Furthermore, we can rearrange the summation to obtain the following result.
\[
  R = \sum_{\ell=0}^{k} \sum_{\nu=0}^{n} \binom{n-\nu}{n-k} 
  \binom{\nu}{\ell} \AP_\ell(x) \, y^{k-\ell}.
\]
Since the terms $\AP_\ell(x) \, y^{k-\ell}$ are independent of the index $\nu$,
we can evaluate the inner sum separately. Using the convolution identity 
\eqref{eq:ud-conv}, this yields
\[
  \sum_{\nu=0}^{n} \binom{n-\nu}{n-k} \binom{\nu}{\ell} 
  = \binom{n+1}{n-k+\ell+1} = \binom{n+1}{k-\ell}.
\]
Finally, 
\[
  R = \sum_{\ell=0}^{k} \binom{n+1}{k-\ell} \AP_\ell(x) \, y^{k-\ell}
  = \sum_{\ell=0}^{k} \binom{n+1}{\ell} \AP_{k-\ell}(x) \, y^\ell,
\]
which is identical to the left-hand side of \eqref{eq:ap-loc}, and we are 
done in this case.

It is easy to see that \eqref{eq:brr} follows from \eqref{eq:brl} by using the
symmetry of the binomial coefficients and rearranging of indices. 
Setting $n = k$, at last \eqref{eq:brl} turns into \eqref{eq:brs}, 
completing the proof.
\end{proofcref}

Corollaries~\ref{cor:main}~--~\ref{cor:main3} follow directly and need no proofs.

\begin{proofcref}{cor:main4}
The first two identities \eqref{eq:brr4} and \eqref{eq:bri4}
follow from Corollary~\ref{cor:main2} by setting 
$x=-\frac{1}{2}$ and $x=1$ in \eqref{eq:brl3}, respectively.
The last equations in \eqref{eq:brc4} are then given by symmetry and 
\eqref{eq:bri4} via
\[
  \sum_{\nu=0}^{n} \binom{n+\nu}{n} 2^{n-\nu}
  = \sum_{\nu=0}^{n} \binom{n+\nu}{\nu} 2^{n-\nu}
  = \sum_{\nu=0}^{n} \binom{2n-\nu}{n-\nu} 2^\nu
  = \sum_{\nu=0}^{n} \binom{2n+1}{\nu} 
  = 2^{2n}.
\]
This completes the proof.
\end{proofcref}

% ------------------------------------------------------------------------------
\subsection{Fibonacci identities}
% ------------------------------------------------------------------------------

\begin{proofcref}{thm:fib}
Let $n \geq k \geq 0$. We use Corollary~\ref{cor:main2} and rewrite \eqref{eq:brl3} as
\[
  \sum_{\nu=0}^{k} \binom{n+1}{k-\nu} x^{\nu} 
  = \sum_{\nu=0}^{k} \binom{n-\nu}{k-\nu} (x+1)^{\nu}.
\]
Applying Proposition~\ref{prop:ident} to the above identity, we infer that 
\eqref{eq:fib-id1}~--~\eqref{eq:fib-id3} provide \eqref{eq:tf1}~--~\eqref{eq:tf3},
respectively. Finally, substracting \eqref{eq:fib-id3} from \eqref{eq:fib-id4} 
implies \eqref{eq:tf4} by considering 
\[
  -(-1)^\nu \FN_\nu - \FN_\nu = -2 \delta_{2 \mid \nu} \, \FN_\nu. \qedhere
\]
\end{proofcref}

\begin{proofcref}{cor:fib}
Let $n \geq 0$. We first observe (cf.~Hoggatt \cite[Sec.\,10, pp.\,52]{Hoggatt:1969}) that
\[
  \FN_{2n+1} - 1 = \sum_{\nu=0}^{n} \FN_{2\nu} \andq
  (-1)^n \FN_{n-1} - 1 = \sum_{\nu=0}^{n} (-1)^\nu \FN_\nu.
  \label{eq:fib-luc}
\]
Actually, the case $n=0$ has to be checked additionally, since indices start 
at $1$ in most references: $\FN_1 - 1 = \FN_{-1} - 1 = \FN_0 = 0$.
The claimed identities then follow from Corollary~\ref{cor:main2} and \eqref{eq:brs3}
by applying Proposition~\ref{prop:ident} with \eqref{eq:fib-id1} and \eqref{eq:fib-id2},
respectively.
\end{proofcref}

% ------------------------------------------------------------------------------
\subsection{Appell and Bernoulli identities}
% ------------------------------------------------------------------------------

\begin{proofcref}{thm:app}
Let ${n \geq k \geq 0}$. We use Theorem~\ref{thm:main} with ${y=1}$. 
From the reflection relation~\eqref{eq:apr}, it follows that 
${\AP_\nu(x+1) = (-1)^\nu \AP_\nu(-x)}$ for ${\nu \geq 0}$. 
Consequently, \eqref{eq:brl} and \eqref{eq:brs} imply the claimed identities 
\eqref{eq:app-ref1} and \eqref{eq:app-ref2}, respectively. 
The values of $\BN_\nu(0)$ and $\EN_\nu(0)$ are given by \eqref{eq:ap-val}.
\end{proofcref}

\begin{proofcref}{cor:app-bp}
We apply Theorem~\ref{thm:app} with ${\AP_\nu(x) = \BN_\nu(x)}$ and 
$\BN_\nu(0) = \BN_\nu$ for $\nu \geq 0$.
First, we show the identities \eqref{eq:app-bp1} and \eqref{eq:app-bp2}.
Let ${n \geq k \geq 0}$. Note that for $k=0$ these identities vanish. 
More precisely, we have $\binom{n+1}{0} - \binom{n}{0} = 0$ on the left-hand sides, 
while the sum vanishes and $\binom{n-1}{-1}=0$ on the right-hand sides, respectively. 
Now, let $k \geq 1$. We use \eqref{eq:app-ref1} and apply relation~\eqref{eq:bp-sym}.
This implies \eqref{eq:app-bp1} and consequently \eqref{eq:app-bp2} by taking $x=0$.
Second, we show the remaining identities \eqref{eq:app-bp3} and \eqref{eq:app-bp4}.
Let $n \geq 0$. Using \eqref{eq:bp-sym}, identity \eqref{eq:app-ref2} turns into 
\eqref{eq:app-bp3}.
Finally, \eqref{eq:app-bp4} easily follows from \eqref{eq:app-bp3} with ${x=0}$, 
where the sum reduces to $\delta_{n > 0}$ on the right-hand side.
This completes the proof. 
\end{proofcref}

\begin{proofcref}{cor:app-bn}
First, we show \eqref{eq:app-bn1}. Let $n \geq k \geq 3$, where $k$ is odd. 
Since the numbers $\BN_\nu$ vanish for odd ${\nu \geq 3}$, we rewrite 
\eqref{eq:app-bp2} of Corollary~\ref{cor:app-bp}, as follows. 
Using the substitution 
\[
  \beta_\nu = \left( \binom{n+1}{k-\nu} - \binom{n-\nu}{k-\nu} \right) \BN_\nu,
\]
we then obtain
\begin{equation} \label{eq:app-bt}
  \beta_{k-1} + \sum_{\substack{\nu=2\\2 \mids \nu}}^{k-3} \beta_\nu 
  = \binom{n-1}{k-1} - \beta_0 - \beta_1. 
\end{equation}
Note that the above sum vanishes for $k=3$, but the equation remains valid in 
that case. Evaluating terms, we derive that
\[
  \beta_{k-1} = \left( \binom{n+1}{1} - \binom{n-k+1}{1} \right) \BN_{k-1} 
  = k \, \BN_{k-1}
\]
and using Table~\ref{tab:ap-val} that
\[
  \beta_0 = \binom{n+1}{k} - \binom{n}{k} = \binom{n}{k-1} \andq  
  \beta_1 = -\frac{1}{2} \left( \binom{n+1}{k-1} - \binom{n-1}{k-1} \right).
\]
The right-hand side of \eqref{eq:app-bt} then becomes
\[
  \frac{1}{2} \binom{n+1}{k-1} - \binom{n}{k-1} + \frac{1}{2} \binom{n-1}{k-1}
  = \frac{1}{2} \Delta^{\!2} \binom{n-1}{k-1} = \frac{1}{2} \binom{n-1}{k-3}.
\]
Putting all together, we infer that \eqref{eq:app-bt}, multiplied by $2$, 
implies \eqref{eq:app-bn1}.

Second, we show \eqref{eq:app-bn2}. Let $n \geq 0$. We use \eqref{eq:app-bp4}
in modified form, since the coefficient of $\BN_n$ vanishes unfortunately.
Therefore, we simply add \eqref{eq:bn-rec} to \eqref{eq:app-bp4}. This yields
\[
  \sum_{\nu=0}^{n} \left( \binom{n+2}{\nu+1} - 1 \right) \BN_\nu = 1.
\]
Now, let $n \geq 2$ be even. Similar to \eqref{eq:app-bt}, 
we rewrite the above equation as
\begin{equation} \label{eq:app-bt2}
  (n+1) \BN_n + \sum_{\substack{\nu=2\\2 \mids \nu}}^{n-2} 
  \left( \binom{n+2}{\nu+1} - 1 \right) \BN_\nu = R,
\end{equation}
where
\[
  R = 1 - \left( \binom{n+2}{1} - 1 \right) + \frac{1}{2} \left( \binom{n+2}{2} - 1 \right)
  = \frac{1}{2} \binom{n}{2}.
\]
Finally, \eqref{eq:app-bt2}, multiplied by $2$, provides \eqref{eq:app-bn2},
finishing the proof.
\end{proofcref}

% ------------------------------------------------------------------------------
\subsection{Euler polynomials and Genocchi numbers}
% ------------------------------------------------------------------------------

\begin{proofcref}{cor:app-ep}
We apply Theorem~\ref{thm:app} with ${\AP_\nu(x) = \EN_\nu(x)}$ and 
$\EN_\nu(0) = \GD_{\nu+1}$ for~${\nu \geq 0}$.
Let ${n \geq k \geq 0}$. Using \eqref{eq:app-ref1} and relation~\eqref{eq:ep-sym},
we get \eqref{eq:app-ep1}, and \eqref{eq:app-ep2} by setting $x=0$.
Finally, the identities \eqref{eq:app-ep3} and \eqref{eq:app-ep4} then follow 
from setting $n=k$, respectively.
\end{proofcref}

\begin{proofcref}{cor:app-gn}
First, we show \eqref{eq:app-gn1}. Let $n \geq k \geq 1$, where $k$ is odd. 
As ${\GD_\nu = 0}$ for odd ${\nu \geq 3}$, we rewrite \eqref{eq:app-ep2} of 
Corollary~\ref{cor:app-ep}. Using the substitution
\[
  \gamma_\nu = \left( \binom{n+1}{k-\nu} + \binom{n-\nu}{k-\nu} \right) \GD_{\nu+1},
\]
we infer from Table~\ref{tab:ap-val} that
\[
  \gamma_k + \sum_{\substack{\nu=1\\2 \mids \nu}}^{k-2} \gamma_\nu 
  = 2 \binom{n}{k} - \gamma_0 = \binom{n}{k} - \binom{n+1}{k} = - \binom{n}{k-1} 
\]
with $\gamma_k = 2 \GD_{k+1}$, implying \eqref{eq:app-gn1}.
Finally, for $n=k$ and shifting the indices $n+1 \mapsto n$ and $\nu+1 \mapsto \nu$ 
in \eqref{eq:app-gn1} gives \eqref{eq:app-gn2}, which is valid for even $n \geq 2$.  
This completes the proof.
\end{proofcref}

% ------------------------------------------------------------------------------
\subsection{Genocchi numbers again}
% ------------------------------------------------------------------------------

The results in this subsection are based on applying the derivative 
$\AP{\!'}_\nu(x)$ in Theorem~\ref{thm:app}.

\begin{proofcref}{cor:app-dgn}
Let $n \geq k \geq 0$. We use Theorem~\ref{thm:app}. Taking the derivative of 
\eqref{eq:app-ref1} and using \eqref{eq:apd} provide
\[
  \sum_{\nu=1}^{k} \binom{n+1}{k-\nu} \nu \AP_{\nu-1}(x) 
  = \sum_{\nu=1}^{k} \binom{n-\nu}{k-\nu} (-1)^{\nu+1} \nu \AP_{\nu-1}(-x). 
\]
Setting $\AP_\nu(x) = \EN_\nu(x)$ and $x=0$, this gives $\nu \EN_{\nu-1}(0) = \GN_\nu$
for $\nu \geq 1$. Since $\GN_0 = 0$ by Table~\ref{tab:ap-val}, the above equation 
is compatible with \eqref{eq:app-dgn1}. Setting $n=k$ finally implies \eqref{eq:app-dgn2}. 
\end{proofcref}

\begin{proofcref}{cor:app-dgn2}
Let $n \geq k \geq 0$. We use Corollary~\ref{cor:app-dgn}.
As $\GN_1 = 1$ by Table~\ref{tab:ap-val} and $\GN_\nu = 0$ for odd $\nu \geq 3$, 
we infer from \eqref{eq:app-dgn1} that
\[
  \sum_{\nu=0}^{k} \binom{n+1}{k-\nu} \GN_\nu
  = \sum_{\nu=0}^{k} \binom{n-\nu}{k-\nu} (-\GN_\nu) \;+\; \ST_k,
\]
where
\[
  \ST_k = 2 \binom{n-1}{k-1} \GN_1 = 2\binom{n-1}{k-1}.
\]
By definition, $\ST_k = 0$ for $k=0$ as required, since $\GN_0 = 0$. 
This shows \eqref{eq:app-dgn3}. For $n=k$, we get \eqref{eq:app-dgn4} with 
$\ST_n = 2\delta_{n > 0}$, completing the proof.
\end{proofcref}

\begin{proofcref}{cor:app-dgn3}
Let $n \geq k \geq 2$ and $k$ be even. We use Corollary~\ref{cor:app-dgn2}
and rewrite \eqref{eq:app-dgn3} with the substitution
\[
  \gamma_\nu = \left( \binom{n+1}{k-\nu} + \binom{n-\nu}{k-\nu} \right) \GN_{\nu}.
\]
By $\GN_0 = 0$, $\GN_1 = 1$, and $\GN_\nu = 0$ for odd $\nu \geq 3$, we obtain
\[
  \gamma_k + \sum_{\substack{\nu=2\\2 \mids \nu}}^{k-2} \gamma_\nu 
  = 2\binom{n-1}{k-1} - \gamma_1 = \binom{n-1}{k-1} - \binom{n+1}{k-1}, 
\]
where $\gamma_k = 2\GN_k$. This implies \eqref{eq:app-dgn5}.
Finally, setting $n=k$ yields \eqref{eq:app-dgn6}, finishing the proof.
\end{proofcref}

% ------------------------------------------------------------------------------
\subsection{Mixed identities}
% ------------------------------------------------------------------------------

Remember the definition of the reciprocal Appell polynomials by \eqref{eq:ap-rec-def}.

\begin{proofcref}{thm:mix} 
Let $n \geq k \geq 0$ and $s \in \CC \setminus \set{0}$. We use Theorem~\ref{thm:main}. 
First, we set $x=0$ and $y=1/s$ in \eqref{eq:brl} and obtain
\begin{align*}
  \sum_{\nu=0}^{k} \binom{n+1}{k-\nu} s^{\nu-k} \AP_\nu(0)
  &= \sum_{\nu=0}^{k} \binom{n-\nu}{k-\nu} s^{\nu-k} \AP_\nu(s^{-1}). 
\shortintertext{Multiplying by $s^k$ and using \eqref{eq:ap-rec-def} yield}
  \sum_{\nu=0}^{k} \binom{n+1}{k-\nu} s^\nu \AP_\nu(0)
  &= \sum_{\nu=0}^{k} \binom{n-\nu}{k-\nu} \AR_\nu(s). 
\end{align*}
The above equation is also valid for $s=0$. 
Replacing $s$ by $x$ provides \eqref{eq:mix-ap1}.

Second, we set $x=1/s$ and $y=-1/s$ in \eqref{eq:brl}. This gives
\begin{align*}
  \sum_{\nu=0}^{k} \binom{n+1}{k-\nu} (-s)^{\nu-k} \AP_\nu(s^{-1})
  &= \sum_{\nu=0}^{k} \binom{n-\nu}{k-\nu} (-s)^{\nu-k} \AP_\nu(0). 
\shortintertext{Multiplying by $(-s)^k$ and using \eqref{eq:ap-rec-def} yield}
  \sum_{\nu=0}^{k} \binom{n+1}{k-\nu} (-1)^\nu \AR_\nu(s)
  &= \sum_{\nu=0}^{k} \binom{n-\nu}{k-\nu} (-s)^\nu \AP_\nu(0). 
\end{align*}
Again, the above equation is also valid for $s=0$, and replacing $s$ by $x$ 
shows \eqref{eq:mix-ap2}. The remaining identities \eqref{eq:mix-ap3} and
\eqref{eq:mix-ap4} follow from setting $n=k$, respectively. 
This completes the proof.
\end{proofcref}

\begin{proofcref}{cor:mix-hs} 
We set $x=0$ in Theorem~\ref{thm:mix}. By \eqref{eq:ap-rec0}, we can replace 
terms by $\AP_0(0) = \AR_\nu(0)$ in \eqref{eq:mix-ap1} and \eqref{eq:mix-ap2}.
By assumption that $\AP_0(0) \neq 0$, we can divide both equations by $\AP_0(0)$, 
yielding the ordinary hockey-stick identities \eqref{eq:hsl} and \eqref{eq:asl}, 
respectively.
\end{proofcref}

\begin{proofcref}{cor:mix}
Let $n \geq k \geq 0$. We apply Theorem~\ref{thm:mix} with $x=2$. By assumption,
$(\AP_\nu(x))_{\nu \geq 0}$ satisfies the reflection relation~\eqref{eq:apr}.
By \eqref{eq:ap-rec1}, $\AR_n(2)$ then vanishes for odd $n \geq 1$. 
As a result, \eqref{eq:mix-ap1} and \eqref{eq:mix-ap2} become \eqref{eq:mix2-ap1} 
and \eqref{eq:mix2-ap2}, as well as, \eqref{eq:mix2-ap4} and \eqref{eq:mix2-ap5} 
by setting $n=k$, respectively. The special values of $\BN_\nu(x)$ and $\EN_\nu(x)$ 
are given by \eqref{eq:ap-val} and \eqref{eq:ap-rec2}. It remains to show
the identities \eqref{eq:mix2-ap3} and \eqref{eq:mix2-ap6}, and to compute the
values $\ST_k$ and $\ST_n$ using Table~\ref{tab:ap-val}, respectively.

Case (1): Regarding \eqref{eq:mix2-ap2} and \eqref{eq:mix2-ap3}, we have to
consider the terms $(-2)^\nu \, \BN_\nu$.
As $\BN_1 = -\frac{1}{2}$ and $\BN_\nu = 0$ for odd $\nu \geq 3$, 
we infer that $\alpha = 0$ and
\[
  \ST_k = -2 \binom{n-1}{k-1} 2 \BN_1 = 2\binom{n-1}{k-1},
\]
which is compatible with $\ST_k = 0$ when $k=0$. 
This implies $\ST_n = 2\delta_{n > 0}$ in \eqref{eq:mix2-ap6}.

Case (2): Similarly as above, we have to consider the terms $(-2)^\nu \GD_{\nu+1}$.
Since $\GN_1 = 1$ and $\GN_\nu = 0$ for odd $\nu \geq 3$, 
we derive that $\alpha = 1$, 
\[
  \ST_k = 2 \binom{n}{k} \GN_1 = 2\binom{n}{k}, \andq \ST_n = 2.
\]
This completes the proof.
\end{proofcref}

\begin{proofcref}{cor:mix-en}
Let $n \geq 0$. From Corollary~\ref{cor:mix}(2) and \eqref{eq:mix2-ap4}, 
it follows that
\[
  \sum_{\nu=0}^{n} \EN_\nu 
  = \sum_{\nu=0}^{n} \binom{n+1}{\nu+1} 2^\nu \GD_{\nu+1} 
  = \frac{1}{2} \sum_{\nu=1}^{n+1} \binom{n+1}{\nu} 2^\nu \GD_{\nu},
\]
showing \eqref{eq:en-sum}. Applying the $\Delta$ operator to the above equation 
for $n-1$ implies that
\[
  \EN_n = \frac{1}{2} \sum_{\nu=1}^{n+1} \binom{n}{\nu-1} 2^\nu \GD_{\nu}
  = \frac{1}{2(n+1)} \sum_{\nu=1}^{n+1} \binom{n+1}{\nu} 2^\nu \GN_{\nu},
\]
giving \eqref{eq:en-sum2}.
\end{proofcref}

\begin{remark}
Note that \eqref{eq:en-sum2} can be derived elementary from the translation 
formula \eqref{eq:apt} for $\EN_n(x+y)$ with $x=0$ and $y = \tfrac{1}{2}$,
and using \eqref{eq:ap-rec2} for $\EN_n$. Hence,
\[
  \EN_n = \sum_{\nu=0}^{n} \binom{n}{\nu} 2^\nu \, \EN_\nu(0)
  = \sum_{\nu=0}^{n} \binom{n}{\nu} 2^\nu \GD_{\nu+1}.
\]
\end{remark}

% ------------------------------------------------------------------------------
\subsection{Power sum identities}
% ------------------------------------------------------------------------------

\begin{proofcref}{thm:sum}
We apply Theorem~\ref{thm:main} with ${\AP_\nu(x) = \BN_\nu(x)}$ for ${\nu \geq 0}$. 
By \eqref{eq:sn-int}, integrating \eqref{eq:brl} over the interval $[x_1,x_2]$ yields
\[
  \int_{x_1}^{x_2} \! \BN_n(t) \, dt = S_\nu(x_1,x_2) \andq
  \int_{x_1}^{x_2} \! \BN_n(t+y) \, dt = \int_{x_1+y}^{x_2+y} \! \BN_n(t) \, dt 
  = S_\nu(x_1+y,x_2+y),
\]
showing the first result. From setting $n=k$, the second result follows. 
\end{proofcref}

\begin{proofcref}{cor:sum}
Using Theorem~\ref{thm:sum} with ${x_1=0}$, ${x_2=x}$, and ${y=1}$, we then 
obtain for ${\nu \geq 0}$ that $S_\nu(0,x) = S_\nu(x)$ and 
$S_\nu(1,x+1) = S_\nu(x) + x^\nu - S_\nu(1)$ by \eqref{eq:sn-delta}. 
Since $S_\nu(1) = 0^\nu$, the term $x^\nu - S_\nu(1)$ vanishes for $\nu=0$, 
implying the first identity. The second identity is finally given by setting $n=k$.  
\end{proofcref}

% ==============================================================================
\section{Conclusion}
% ==============================================================================

The present paper was motivated by a similar identity 
(see \cite[Thm.\,2.7]{Kellner:2023}) of the author that reads as follows.
For $n, r \geq 1$, we have
\[
  \sum_{k=0}^{r-1} \binom{n+r}{k} (-1)^k x^{r-1-k} 
  = \sum_{k=0}^{r-1} \binom{n+k}{k} (-1)^k (x-1)^{r-1-k}.
\]
The above identity occurred in the context of certain binomial sums that are 
conjectured to be nonintegral. The search for the origin of this identity eventually 
led to the generalized identities with the Appell polynomials in this paper. 
From Proposition~\ref{prop:ident}, we casually derive the following connection 
with the Fibonacci numbers.

\begin{corollary}
For $n, r \geq 1$, we have
\begin{align*}
  \sum_{k=0}^{r-1} \binom{n+r}{k} (-1)^k \FN_{2(r-1-k)} 
  &= \sum_{k=0}^{r-1} \binom{n+k}{k} (-1)^k \FN_{r-1-k}, \\
  \sum_{k=0}^{r-1} \binom{n+r}{k} (-1)^{r-k} \FN_{r-1-k} 
  &= \sum_{k=0}^{r-1} \binom{n+k}{k} \FN_{r-1-k}, \\
  \sum_{k=0}^{r-1} \binom{n+r}{k} \FN_{r-1-k} 
  &= \sum_{k=0}^{r-1} \binom{n+k}{k} \FN_{2(r-1-k)}.
\end{align*}
\end{corollary}

Finally, it turns out that the identity (see \cite[Sec.\,5.1, Eq.\,(5.19), p.\,166]{GKP:1994})
\[
  \sum_{k=0}^{m} \binom{m+r}{k} x^k y^{m-k} 
  = \sum_{k=0}^{m} \binom{-r}{k} (-x)^k (x+y)^{m-k}
\]
can be shown to be equivalent to \eqref{eq:brl2} of Corollary~\ref{cor:main}.
However, the identity above was called a \emph{curious relationship} 
by the authors of the latter reference.

In a sequel paper \cite{Kellner:2026}, the boomerang identities are considered
in the context of sequences in rings and their binomial transforms.

% ==============================================================================
\appendix
\section{Summary}
% ==============================================================================

Let $n \geq k \geq 0$ and $x, y \in \CC$. Recalling the notation of Section~\ref{sec:appl}, 
we summarize the obtained boomerang identities of the paper by
\begin{align*}
  \sum_{\nu=0}^{k} \binom{n+1}{k-\nu} \, \RT_\nu
  &= \sum_{\nu=0}^{k} \binom{n-\nu}{k-\nu} \, \CT_\nu \;+\; \ST_k 
\shortintertext{and}  
  \sum_{\nu=0}^{n} \binom{n+1}{\nu+1} \, \RT_\nu
  &= \sum_{\nu=0}^{n} \CT_\nu \;+\; \ST_n, 
\end{align*}
where $\RT$, $\CT$, and $\ST$ represent terms in the rows and columns of 
Pascal's triangle, as well as, supplementary terms, respectively.

\begin{table}[H] \small
\setstretch{1.22}
\begin{center}
\begin{tabular}{ccc}
  \toprule
  $\RT_\nu$ & $\CT_\nu$ & Application \\
  \midrule
  $\AP_\nu(x) \, y^{k-\nu}$ & $\AP_\nu(x+y) \, y^{k-\nu}$ & generic \\
  $x^{\nu} y^{k-\nu}$ & $(x+y)^{\nu} y^{k-\nu}$ & \ditto \\
  $x^{\nu}$ & $(x+1)^{\nu}$  & \ditto \\
  \midrule
  $0^\nu$ & $1$ & Hockey stick \\
  $(-1)^\nu$ & $0^\nu$ & \ditto \\
  \midrule
  $1$ & $2^\nu$ & sums \\
  $(-\tfrac{1}{2})^\nu$ & $(\tfrac{1}{2})^\nu$ & \ditto \\
  \midrule
  $\FN_\nu$ & $\FN_{2\nu}$ & Fibonacci number \\
  $(-1)^\nu \FN_{2\nu}$ & $(-1)^\nu \FN_\nu$ & \ditto \\
  $(-1)^\nu \FN_\nu$ & $-\FN_\nu$ & \ditto \\
  $(-1)^\nu (\FN_{2\nu} - 2 \delta_{2 \mid \nu} \, \FN_\nu)$ 
  & $-(\FN_{2\nu} - 2 \delta_{2 \mid \nu} \, \FN_\nu)$ & \ditto \\
  \bottomrule
\end{tabular}
\caption{Basic terms.}
\label{tab:term}
\end{center}
\end{table}

\begin{table}[H] \small
\setstretch{1.3}
\begin{center}
\begin{tabular}{ccccc}
  \toprule
  $\RT_\nu$ & $\CT_\nu$ & $\ST_k$ & $\ST_n$ & Application \\
  \midrule
  $\AP_\nu(x)$ & $(-1)^\nu \AP_\nu(-x)$ &&& Appell with reflection \\
  $\BN_\nu(x)$ & $\BN_\nu(x)$ & $\sum_{\nu=1}^{k} \binom{n-\nu}{k-\nu} \nu x^{\nu-1}$ 
  & $\sum_{\nu=1}^{n} \nu x^{\nu-1}$ & Bernoulli polynomial \\
  $\EN_\nu(x)$ & $-\EN_\nu(x)$ & $2 \sum_{\nu=0}^{k} \binom{n-\nu}{k-\nu} x^\nu$ 
  & $2 \sum_{\nu=0}^{n} x^\nu$ & Euler polynomial \\
  \midrule
  $\BN_\nu$ & $(-1)^\nu \BN_\nu$ &&& Bernoulli number \\
  $\BN_\nu$ & $\BN_\nu$ & $\binom{n-1}{k-1}$ & $\delta_{n > 0}$ & \ditto \\
  \midrule
  $\GN_\nu$ & $(-1)^{\nu+1} \GN_\nu$ &&& Genocchi number \\
  $\GN_\nu$ & $-\GN_\nu$ & $2\binom{n-1}{k-1}$ & $2\delta_{n > 0}$ & \ditto \\
  $\GD_{\nu+1} $ & $(-1)^\nu \GD_{\nu+1}$ &&& \ditto \\
  $\GD_{\nu+1} $ & $-\GD_{\nu+1}$ & $2\binom{n}{k}$ & $2$ & \ditto \\
  \midrule
  $S_\nu(x_1,x_2) \, y^{k-\nu}$ & $S_\nu(x_1+y,x_2+y) \, y^{k-\nu}$ &&& power sum \\
  $S_\nu(x)$ & $S_\nu(x)$ & $\sum_{\nu=1}^{k} \binom{n-\nu}{k-\nu} x^\nu$ 
  & $\sum_{\nu=1}^{n} x^\nu$ & \ditto \\
  \bottomrule
\end{tabular}
\caption{Terms derived from Appell polynomials.}
\label{tab:term2}
\end{center}
\end{table}

\begin{table}[H] \small
\setstretch{1.3}
\begin{center}
\begin{tabular}{ccccc}
  \toprule
  $\RT_\nu$ & $\CT_\nu$ & $\ST_k$ & $\ST_n$ & Application \\
  \midrule
  $x^\nu \! \AP_\nu(0)$ & $\AR_\nu(x)$ &&& Appell polynomial \\
  $(-1)^\nu \AR_\nu(x)$ & $(-x)^\nu \! \AP_\nu(0)$ &&& \ditto \\
  $\AR_\nu(2)$ & $(-2)^\nu \AP_\nu(0)$ &&& Appell with reflection \\
  \midrule
  $2^\nu \BN_\nu$ & $\BN_\nu + \tfrac{1}{2} \GN_\nu$ &&& mixed \\
  $\BN_\nu + \tfrac{1}{2} \GN_\nu$ & $(-2)^\nu \BN_\nu$ &&& \ditto \\
  $\BN_\nu + \tfrac{1}{2} \GN_\nu$ & $2^\nu \BN_\nu$ & $2\binom{n-1}{k-1}$ 
  & $2\delta_{n > 0}$ & \ditto \\
  $2^\nu \GD_{\nu+1}$ & $\EN_\nu$ &&& \ditto \\
  $\EN_\nu$ & $(-2)^\nu \GD_{\nu+1}$ &&& \ditto \\
  $\EN_\nu$ & $-2^\nu \GD_{\nu+1}$ & $2\binom{n}{k}$ & $2$ & \ditto \\
  \bottomrule
\end{tabular}
\caption{Terms derived from reciprocal Appell polynomials.}
\label{tab:term3}
\end{center}
\end{table}

% ==============================================================================
\bibliographystyle{amsplain}

\end{document}